\documentclass[12pt]{article}

\usepackage{amsmath, amsfonts, amssymb}
\usepackage{graphicx}
\usepackage{epstopdf}
\usepackage{bm}
\usepackage{color}
\usepackage{lscape}

\usepackage{algorithm}
\usepackage{algorithmicx}

\newcommand{\bd}[1]{\mathbf{#1}}

\newcommand{\del}{\delta}

\newcommand{\pa}{\partial}
\newcommand{\beq}{\begin{equation}}
\newcommand{\eeq}{\end{equation}}

\DeclareMathOperator\erf{erf}

\title{A Novel Regularization for Higher Accuracy in the Solution of 3D Stokes Flow}
\author{J. Thomas Beale\thanks{Department of Mathematics, Duke University, Durham, NC, 27708 USA beale@math.duke.edu} \and Christina Jones\thanks{Department of Mathematics, Farmingdale State College, SUNY, Farmingdale, NY 11735, USA jonecj2@farmingdale.edu} \and Jillian Reale\thanks{Department of Mathematics, Farmingdale State College, SUNY, Farmingdale, NY 11735, USA realj@farmingdale.edu} \and Svetlana Tlupova\thanks{Department of Mathematics, Farmingdale State College, SUNY, Farmingdale, NY 11735, USA tlupovs@farmingdale.edu}}

\date{\today}

\begin{document}

\maketitle

\begin{abstract} 
Many problems in fluid dynamics are effectively modeled as Stokes flows - slow, viscous flows where the Reynolds number is small. Boundary integral equations are often used to solve these problems, where the fundamental solutions for the fluid velocity are the Stokeslet and stresslet. One of the main challenges in evaluating the boundary integrals is that the kernels become singular on the surface. A regularization method that eliminates the singularities and reduces the numerical error through correction terms for both the Stokeslet and stresslet integrals was developed in Tlupova and Beale, JCP (2019). In this work we build on the previously developed method to introduce a new stresslet regularization that is simpler and results in higher accuracy when evaluated on the surface. Our regularization replaces a seventh-degree polynomial that results from an equation with two conditions and two unknowns with a fifth-degree polynomial that results from an equation with one condition and one unknown. Numerical experiments demonstrate that the new regularization retains the same order of convergence as the regularization developed by Tlupova and Beale but shows a decreased magnitude of the error.
\end{abstract}


\section{Introduction}

Many problems in fluid dynamics are modeled as Stokes flows - particle interactions in slow, viscous flow that results in a small Reynolds number. The equations that describe these flows, the incompressible Stokes equations, are
\beq
	-\nabla p + \Delta \bd{u} = 0,\qquad \nabla\cdot \bd{u} = 0,
\eeq
where $p$ is the fluid pressure and $\bd{u}$ is the fluid velocity. The Stokeslet and stresslet are the primary fundamental solutions for the fluid velocity,
\begin{subequations}
\begin{align}
	\label{stokeslet}
	S_{ij}(\bd{y,x}) &= \frac{\del_{ij}}{|\bd{y} - \bd{x}|} + \frac{(y_i - x_i)(y_j - x_j)}{|\bd{y} - \bd{x}|^3}, \\[6pt]
	\label{stresslet}
	T_{ijk} (\bd{y,x}) &= -\frac{6(y_i - x_i)(y_j - x_j)(y_k - x_k)}{|\bd{y} - \bd{x}|^5},
\end{align}
\end{subequations}
where $\del_{ij}$ is the Kronecker delta and $i,j,k = 1,2,3$ are Cartesian coordinates, $\bd{x}$ is a source point, and $\bd{y}$ is a target point. When used in boundary integral methods, these lead to the single and double layer representations of Stokes flow, respectively,
\begin{subequations}
\begin{align}
	\label{SingleLayer}
	u_i(\bd{y}) &= \frac{1}{8\pi}\int_{\pa\Omega} S_{ij}(\bd{y,x}) f_j(\bd{x})dS(\bd{x}), \\[6pt]
	\label{DoubleLayer}
	w_i(\bd{y}) &= \frac{1}{8\pi}\int_{\pa\Omega} T_{ijk} (\bd{y,x}) q_j(\bd{x}) n_k(\bd{x})dS(\bd{x}),
\end{align}
\end{subequations}
where $n_k$ are the components of the unit outward normal vector to the surface $\pa\Omega$ of a bounded domain $\Omega$. The integral in \eqref{SingleLayer} is continuous across $\pa\Omega$, while the integral in \eqref{DoubleLayer} is discontinuous and has a jump of $\mp 4\pi q_i(\bd{x}_0)$ in the limit from either the interior or exterior of the domain.

Computing the layer representations~\eqref{SingleLayer}-\eqref{DoubleLayer} requires addressing the singularities that develop as $r = |\bd{x}-\bd{y}|$ approaches zero. When evaluating the integrals near the surface, the kernels are nearly singular and straightforward quadratures fail to capture them accurately. Tlupova and Beale \cite{tlupova-beale} introduced regularizations for the Stokeslet and stresslet that result in high accuracy when evaluating at points on and off the surface. The method is based on smoothing the kernels using a regularization parameter $\del>0$ developed for the Laplace kernels in~\cite{beale04, beale-lai}, then applying a simple quadrature  of \cite{wilson, beale16}. For the nearly singular case, corrections are added to reduce the regularization error to the $O(\del^3)$ terms.

When evaluating the integrals on the boundary, e.g. when solving integral equations, special smoothing functions are designed \cite{tlupova-beale} that achieve $O(\del^5)$ accuracy without requiring corrections. In addition, these regularizations do not require that adjustments be made to the grid around the singularity.

In this paper we introduce a new smoothing function for the stresslet~\eqref{DoubleLayer} that results in higher accuracy in the computation of the stresslet at points on the surface. We first summarize the method of \cite{tlupova-beale} in Section 2. The new regularization for the stresslet is developed in Section 3. The results of numerical experiments using the original and new regularizations for three surfaces - a sphere, an ellipsoid, and a four-atom molecular surface - are presented in Section 4.

\section{Numerical Method}
\label{sec:numerical_method}

We now briefly describe the idea of regularization from \cite{tlupova-beale}. We demonstrate the main concepts on the stresslet~\eqref{DoubleLayer} as this is the focus of this paper; the approach for the Stokeslet~\eqref{SingleLayer} is similar, for details we refer the reader to~\cite{tlupova-beale}. First, the singularity is reduced in the stresslet through subtraction, resulting in,
\begin{equation}
	\label{DoubleLayer_subtract}
	w_i(\bd{y}) = \frac{1}{8\pi}\int_{\pa\Omega} T_{ijk} (\bd{y,x}) [q_j(\bd{x}) - q_j(\bd{x}_0)] n_k(\bd{x}) dS(\bd{x}) + \frac{1}{8\pi}\chi (\bd{y}) q_i(\bd{x}_0),
\end{equation}
where $\bd{x}_0$ is the boundary point closest to $\bd{y}$, and we have applied the well known identity (see, for example, \cite{pozrikidis92} sec. 2.1-2.3)
\begin{equation}
	\int_{\pa\Omega} T_{ijk} (\bd{y,x}) n_k(\bd{x})dS(\bd{x}) = \chi (\bd{y})\del_{ij},
\end{equation}
where $\chi (\bd{y}) = 8\pi, 4\pi, 0$ if $\bd{y}$ is inside, on, and outside the boundary, respectively. 

The stresslet is then regularized 
\begin{align}
	\label{StressletHH}
	\bd{w}_\del(\bd{y}) = -\frac{3}{4\pi} \int_{\pa\Omega} &[(\bd{y}-\bd{x})\cdot \tilde{\bd{q}}(\bd{x})] [(\bd{y}-\bd{x})\cdot \bd{n}(\bd{x})] (\bd{y}-\bd{x}) \frac{s_3(r/\del)}{r^5} dS(\bd{x})  + \frac{1}{8\pi}\chi (\bd{y}) \bd{q}(\bd{x}_0), 
\end{align}
where $\tilde{\bd{q}}(\bd{x}) = \bd{q}(\bd{x})-\bd{q}(\bd{x}_0)$, and  $s_3$ is chosen with $\lim_{\rho\to\infty}s_3(\rho)=1$, $s_3(\rho)=O(\rho^5)$ for small $\rho$, and $s_3(r/\del)/r^5$ smooth for a fixed parameter $\delta>0$. Once the integrands are smoothed out, the integrals are discretized using the quadrature method for closed surfaces introduced in \cite{wilson} and explained in \cite{beale16}. The error due to regularization is $O(\del)$, and correction terms were derived analytically in~\cite{tlupova-beale} to reduce the $O(\del)$ and $O(\del^2)$ terms, resulting in the final computation accurate to $O(\del^3)$. 

For the case of solving the stresslet at points on the surface, such as when solving integral equations, a special regularization can be designed to achieve high accuracy to $O(\del^5)$ without the need to compute corrections. In~\cite{tlupova-beale}, such a smoothing function was found by setting $s^\sharp_3 (r) = s_3(r) + a r s'_3(r) + b r^2 s{''}_3(r)$, with $a$ and $b$ being constants chosen to make two moments involving $s_3$ equal to 0. The resulting smoothing function for solving the stresslet on the surface given in~\cite{tlupova-beale} is
\beq
	\label{s3_high}
	s_3^{\#}(r) = \erf(r) - \frac{2}{9} r (9 + 6 r^2 - 36 r^4  + 8 r^6) e^{-r^2}/\sqrt{\pi},
\eeq
where $\erf$ is the error function.

As discussed in~\cite{tlupova-beale}, the error in the double layer integral evaluated on the surface using the smoothing~\eqref{s3_high} is expected to behave as
\beq 
	\label{error_estimate}
	\epsilon_w \leq C_1\delta^5 + C_2 h^2\,e^{-c_0(\delta/h)^2},
\eeq
where $h$ is the grid spacing chosen in coordinate planes for the discretization of the integrals. The first term is due to regularizing the kernels, and the second term is due to discretizing the integrals. As such, the accuracy depends critically on the relationship between $\delta$ and $h$. A large enough choice of $\delta$ is needed to ensure the regularization error is dominant over the discretization error, so that the total error approaches $O(h^5)$. We generally take $\delta/h =$ constant for simplicity and in practice, $\delta/h=3$ works well to maintain the high order in the regularization error.

\section{New Regularization}
\label{sec:new_regularization}

The new regularization we propose increases the accuracy of evaluating the stresslet on the surface by using a slightly simpler smoothing function in place of \eqref{s3_high}. As mentioned earlier, the special smoothing was found in~\cite{tlupova-beale} by setting two moment conditions to 0. We have determined however, that one condition will suffice. Specifically, in the original derivation in~\cite{tlupova-beale}, when evaluating at points on the surface we have $\lambda=0$ thus making the condition requiring (40b) equals to 0 unnecessary; see below. This leaves only one moment condition where a similar integral with $\eta^7$ in place of $\eta^5$ is equal to 0. This allows us to create the new smoothing function by setting $s_3^{\#}(r) = s_3(r) + a r s'_3(r)$, and solving for $a$. 

We start with the original smoothing function from~\cite{tlupova-beale},
\beq
	s_3(r) = \erf(r) - 2 r (\frac{2}{3}r^2 + 1) e^{-r^2}/\sqrt{\pi},
\eeq
and compute 
\beq
	r\, s'_3(r) = \frac{8}{3\sqrt{\pi}}\, r^5\, e^{-r^2}.
\eeq
The integral moment condition is
\beq
	\label{Integral_condition}
	\int_{0}^\infty r^2(s_3^{\#}(r) - 1) dr = 0.
\eeq
Since
\beq
	\int_{0}^\infty r^2(s_3(r) - 1) dr = -\frac{8}{3\sqrt{\pi}}, \qquad \int_{0}^\infty r^2(r s'_3(r)) dr = \frac{8}{\sqrt{\pi}},
\eeq
we can therefore set
\beq
	s_3^{\#} = s_3 + \frac{1}{3}r s'_3
\eeq
to satisfy the integral condition~\eqref{Integral_condition}, which leads to the smoothing function
\beq
	\label{s3_high_new}
	s_3^{\#}(r) = \erf(r) - \frac{2}{9} r (9 + 6 r^2 - 4 r^4) e^{-r^2}/\sqrt{\pi}.
\eeq
Note that the polynomial term in this new function has highest power $r^5$, whereas the original function in~\eqref{s3_high} has $r^7$. 

In the derivation of~\eqref{s3_high} in~\cite{tlupova-beale}, the moment condition~\eqref{Integral_condition} was imposed as well as the zero moment condition
\beq 
	\int_0^\infty (s_3^\sharp(r) - 1)dr = 0.
\eeq
However, for the stresslet integral in the subtracted form~\eqref{DoubleLayer_subtract}, the contribution of this moment to the integral is zero, so that this condition can be omitted. 
More generally, for an integral not in the subtracted form, the original version~\eqref{s3_high} could be used.
We obtain~\eqref{s3_high} in the manner described for~\eqref{s3_high_new}, but with 
$s_3^\sharp$ in the form $s_3 + ars_3' + br^2s_3''$ and $a,b$ chosen
to satisfy the two conditions.
  The situation is analogous to that for the simpler case of the double layer potential for a harmonic function; see p. 607 of~\cite{beale04}.
  

\section{Numerical experiments}
\label{sec:numerical_experiments}

We performed numerical experiments to test the new regularization using three surfaces: a unit sphere, an ellipsoid, and a four-atom molecular surface,

\begin{subequations}
\begin{align}
	\label{Sphere}
		\phi(x_1,x_2,x_3) &=  x_1^2 + x_2^2 + x_3^2 - 1, \\[6pt]
	\label{Ellipsoid}
		\phi(x_1,x_2,x_3) &= \frac{x_1^2}{a^2} + \frac{x_2^2}{b^2} + \frac{x_3^2}{c^2} - 1, \\[6pt]
	\label{Molecule}
		\phi(x_1,x_2,x_3) &= \sum_{k=1}^4 \exp(- |{\bf x} - {\bf x}_k|^2/r^2) - c.
\end{align}
\end{subequations}

For the ellipsoid~\eqref{Ellipsoid} we set $a=1,b=0.6,c=0.4$, and for the molecule surface~\eqref{Molecule} we use centers $\bd{x}_1 = (\sqrt{3}/3,0,-\sqrt{6}/12)$, $\bd{x}_{2,3} = (-\sqrt{3}/6,\pm .5,-\sqrt{6}/12)$, $\bd{x}_4 = (0,0,\sqrt{6}/4)$ and $r = .5$, $c = .6$, as in~\cite{beale16}. The number of quadrature points generated to represent each surface for different grid sizes $h$ are listed in Table~\ref{Quad_pts}.

\begin{table}[htb]
\centering
\begin{tabular}{c|r|r|r}
\hline
$h$ & Sphere & Ellipsoid & Molecule \\
\hline
1/32 & 17070   & 6902  & 9562 \\
\hline
1/64 & 68166   & 27566  & 38354 \\
\hline
1/128 & 272718   & 110250  & 153399 \\
\hline
\end{tabular}
\caption{Number of quadrature points for the unit sphere; ellipsoid $a=1,b=0.6,c=0.4$; and the molecular surface from \cite{beale16}.}
\label{Quad_pts}
\end{table}


\subsection{Sum of single and double layer}

One of the advantages of using boundary integral formulations is that jumps in the physical quantities across interfaces get incorporated into the integrals naturally. Specifically, the general integral formulation, expressed as the sum of the single and double layer integrals,
\beq 
	\label{Sum_SLDL}
	u_i(\bd{y}) = 
	-\frac{1}{8\pi}\int_{\partial\Omega} S_{ij} (\bd{y,x}) [f]_j(\bd{x}) dS(\bd{x}) 
	- \frac{1}{8\pi}\int_{\partial\Omega} T_{ijk} (\bd{y,x}) [u]_j(\bd{x}) n_k(\bd{x})dS(\bd{x}),
\eeq
has $[f] = f^+ - f^- = (\sigma^+ - \sigma^-)\cdot \bd{n} $ as the jump in surface force and $[u]$ as the jump in velocity. Here $\bd{n}$ is the outward unit normal, and the plus/minus signs denote the outside/inside of the boundary. We use the following solution from \cite{tlupova-beale}. On the inside, we assume the velocity is given by a point force singularity of strength $\bd{b} = (1,0,0)$, placed at $\bd{y}_0 = (2,0,0)$. The solution is given by the Stokeslet velocity
\beq 
	\label{Stokeslet_point}
	u^-_i(\bd{y}) = \frac{1}{8\pi} S_{ij}b_j = \frac{1}{8\pi} \Big( \frac{\delta_{ij}}{r} + \frac{\hat{y}_i \hat{y}_j}{r^3} \Big)b_j,
\eeq 
and the stress tensor is 
\beq
	\label{Stress_point}
	\sigma^-_{ik}(\bd{y}) = \frac{1}{8\pi} T_{ijk} b_j = \frac{-6}{8\pi}\frac{\hat{y}_i \hat{y}_j \hat{y}_k}{r^5} b_j,
\eeq
where $\hat{\bd{y}} = \bd{y}-\bd{y}_0$, $r=|\hat{\bd{y}}|$. {We assume this data for the inside of the boundary, and take the solution to be $u^+=0$, $\sigma^+=0$ for the outside.  The jumps $[u]$ and $[f]$ are evaluated at the quadrature points using these inside/outside values. The exact solution on the boundary is the average of outside and inside, or half of the formula for $u_i$ in \eqref{Stokeslet_point}. 

We define the error at a single point as $e(\bd{x}) = \lvert\bd{u}^{\textrm{computed}} (\bd{x}) - \bd{u}^{\textrm{exact}} (\bd{x})\rvert$, where $\lvert\cdot \rvert$ is the vector's Euclidean norm. We then measure either the max or the $L_2$ norm of this error over the evaluation points. The $L_2$ norm is defined as $\lVert e\rVert_2 = \big(\sum_\bd{x} e^2(\bd{x})/n\big)^{1/2}$, where $n$ is the number of evaluation points. Figures~\ref{Sum-del3h} and~\ref{Sum-del1h} compare the errors for the three surfaces using the original regularization~\eqref{s3_high} and the new regularization~\eqref{s3_high_new}. Figure~\ref{Sum-del3h} shows the errors using the larger regularization $\del/h=3$, and Figure~\ref{Sum-del1h} shows the errors using the smaller regularization $\del/h=1$, as the grid size $h$ is refined. Following the error estimate in~\eqref{error_estimate}, when the regularization parameter is chosen large enough, such as $\delta/h=3$, the regularization error is larger than the discretization error, and the overall error is estimated at $O(h^5)$. This is observed with the sphere and the molecular surfaces. The thin ellipsoid does not fit the estimate as well due to the larger curvature and varied spacing, expected to improve with grid refinement. The new regularization function~\eqref{s3_high_new} gives smaller errors in all three cases, most dramatically in the case of the ellipsoid. More precisely, we observe an improvement by approximately a factor of two for the sphere, about a factor of six for the ellipsoid, and about a factor of five for the molecule. For the smaller regularization parameter $\delta/h=1$, observed convergence is $O(h)$ and the new regularization does not make a notable difference, so this regularization regime is not recommended in practice.

\begin{figure}[htb]
\centering
\includegraphics[scale=0.4]{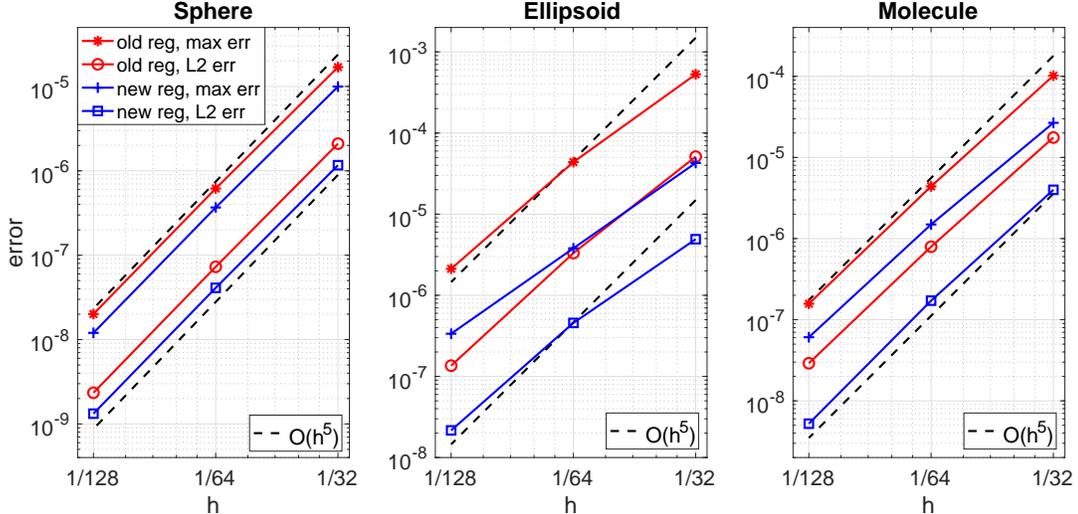}
\caption{Sum of single and double layer, errors over quadrature points for three surfaces; regularization parameter $\del=3h$.}
\label{Sum-del3h}
\end{figure}

\begin{figure}[htb]
\centering
\includegraphics[scale=0.4]{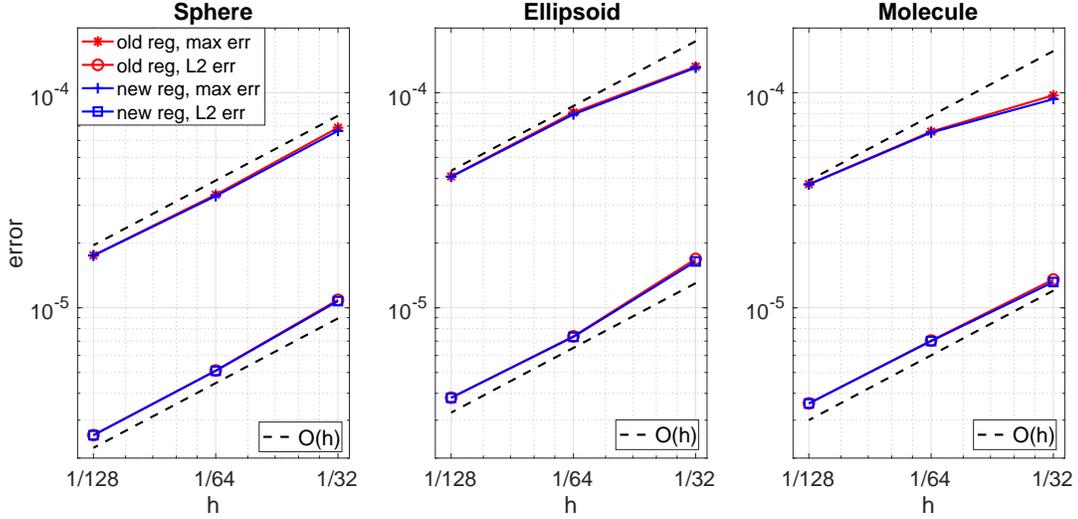}
\caption{Sum of single and double layer, errors over quadrature points for three surfaces; regularization parameter $\del=h$.}
\label{Sum-del1h}
\end{figure}


\subsection{Flow due to an interface with different viscosities}

Here we revisit another example from~\cite{tlupova-beale}, of an interface between two fluids with different viscosities, and an integral equation must be solved to find the interface velocity. The interface undergoes a discontinuity in the surface force $[\bd{f}]$, while the velocity across the interface is continuous \cite{pozrikidis92}. The integral equation for the interface velocity is given by
\begin{align}
	(\lambda+1) u_i(\bd{x}_0) = &-\frac{1}{4\pi \mu_0}\int_{\pa\Omega} S_{ij}(\bd{x}_0,\bd{x}) [f]_j(\bd{x})dS(\bd{x}) \nonumber\\
	&+ \frac{\lambda-1}{4\pi}\int_{\pa\Omega} T_{ijk} (\bd{x}_0,\bd{x}) u_j(\bd{x}) n_k(\bd{x})dS(\bd{x})
\label{IntEq-1surf}
\end{align}
for $\bd{x}_0\in \pa\Omega$, where $\mu_0,\mu_1$ are the external and internal fluid viscosities and $\lambda=\mu_1/\mu_0$. The discontinuity in the surface force is given by $[\bd{f}] = 2\gamma H \bd{n} - \nabla_S\gamma$, where $\gamma$ is the surface tension, $H$ is the mean curvature, and $\bd{n}$ is the outward unit normal \cite{pozrikidis92}. In our numerical tests, we set $\mu_0=1, \mu_1=2$, and $\gamma=1+x_1^2$. We solve the integral equation using successive evaluations, i.e.,
\begin{align}
	(\lambda+1) u^{N}_i(\bd{x}_0) = &-\frac{1}{4\pi \mu_0}\int_{\pa\Omega} S_{ij}(\bd{x}_0,\bd{x}) [f]_j(\bd{x})dS(\bd{x}) \nonumber\\
	&+ \frac{\lambda-1}{4\pi}\int_{\pa\Omega} T_{ijk} (\bd{x}_0,\bd{x}) u^{N-1}_j(\bd{x}) n_k(\bd{x})dS(\bd{x}),
	\label{Successive}
\end{align}
for $N=1,2,...$, and $\bd{u}^0 = \bd{0}$. We stop these iterations when the iteration error, defined as 
\beq
	\label{error_iter}
	e^N := \max_{\bd{x}_0} \lvert \bd{u}^{N}-\bd{u}^{N-1} \rvert,
\eeq 
is below a prescribed tolerance, and $\lvert\cdot \rvert$ is the vector's Euclidean norm. Since the exact solution is not known, we check the convergence rates empirically by defining
\beq
	\label{error_h}
	e_h (\bd{x}) = \bd{u}_{h} (\bd{x}) - \bd{u}_{h/2} (\bd{x}),
\eeq
and taking either the max or the $L_2$ norm of this error over the surface points given by $h$, the larger of the two grid sizes used. These errors are shown in Table~\ref{BIE-1Surf} for the ellipsoid $a=1,b=0.6,c=0.4$, with $\del=3h$. It takes about $N=12$ iterations for the iteration error \eqref{error_iter} to reach below $10^{-10}$. To minimize the error coming from evaluating the single layer integral with the surface tension density (the nonhomogeneous term in~\eqref{Successive}), we compute the single layer integral with increased resolution before solving the integral equation~\eqref{Successive}. Specifically, we solved the integral equation~\eqref{Successive} for each of the values of $h$, but in each case computed the Stokeslet integral at the needed points using the finer grid $h=1/256$. Table~\ref{BIE-1Surf} compares the new regularization with the original one, and shows an improvement when using the new function.

\begin{table}[htb]
\centering
\begin{tabular}{c||c|c||c|c}
\hline
\multicolumn{1}{c||}{} &
\multicolumn{2}{c||}{\ \ \ Original regularization \ \ \ } &
\multicolumn{2}{c}{\ \ \ New regularization \ \ \ } \\
\hline
$h$ & \ \ \ $\lVert e_h\rVert_\infty$ \ \ \ & \ \ \ $\lVert e_h\rVert_2$ \ \ \ & \ \ \ $\lVert e_h\rVert_\infty$ \ \ \ & \ \ \ $\lVert e_h\rVert_2$ \ \ \ \\
\hline
1/16 & 6.93e-03 & 1.85e-03    & 1.69e-03 & 4.77e-04 \\
\hline
1/32 & 6.86e-04 & 1.20e-04    & 2.03e-04 & 3.74e-05 \\
\hline
1/64 & 5.50e-05 & 7.97e-06    & 1.81e-05 & 2.55e-06 \\
\hline
\end{tabular}
\caption{Flow due to an interface, for the ellipsoid $a=1,b=0.6,c=0.4$. The single layer integral computed using $h=1/256$. Grid size $h$, max and $L_2$ norms of the error defined in \eqref{error_h}. Regularization parameter $\del=3h$.}
\label{BIE-1Surf}
\end{table}

\section{Conclusions}
\label{sec:conclusions}

We have introduced a new regularization function for evaluating the double layer potential (stresslet integral) in Stokes flows at points on the surface with high accuracy. The new function only requires one moment condition and has a lower degree polynomial as a result. Numerical tests demonstrate that the new regularization retains the same order of convergence as the regularization developed in prior work but shows a decreased magnitude of the error.


\section*{Acknowledgments}

The work of CJ, JR, and ST was supported by the National Science Foundation grant DMS-2012371. The authors thank an anonymous reviewer for very constructive suggestions.


\end{document}